\documentclass[12pt]{article}
\usepackage{amssymb}
\usepackage{amsfonts}
\usepackage{mitpress}

\newtheorem{theorem}{Theorem}
\newtheorem{acknowledgement}[theorem]{Acknowledgement}

\newtheorem{conjecture}[theorem]{Conjecture}

\newenvironment{proof}[1][Proof of Theorem 1]{\noindent\textbf{#1.} }{\ \rule{0.5em}{0.5em}}
\input{tcilatex}

\begin{document}

\title{Online version of the theorem of Thue}
\author{Jaros\l aw Grytczuk, Piotr Szafruga, Micha\l\ Zmarz \\
Faculty of Mathematics and Computer Science, Jagiellonian University, 30-348
Krak\'{o}w, Poland; grytczuk@tcs.uj.edu.pl, szfruga@tcs.uj.edu.pl,
zmarz@tcs.uj.edu.pl}
\maketitle

\begin{abstract}
A sequence $S$ is \emph{nonrepetitive} if no two adjacent blocks of $S$ are
the same. In 1906 Thue proved that there exist arbitrarily long
nonrepetitive sequences over $3$ symbols. We consider the online variant of
this result in which a nonrepetitive sequence is constructed during a play
between two players: Bob is choosing a position in a sequence and Alice is
inserting a symbol on that position taken from a fixed set $A$. The goal of
Bob is to force Alice to create a repetition, and if he succeeds, then the
game stops. The goal of Alice is naturally to avoid that and thereby to
construct a nonrepetitive sequence of any given length.

We prove that Alice has a strategy to play arbitrarily long provided the
size of the set $A$ is at least $12$. This is the online version of the
Theorem of Thue. The proof is based on nonrepetitive colorings of
outerplanar graphs. On the other hand, one can prove that even over $4$
symbols Alice has no chance to play for too long. The minimum size of the
set of symbols needed for the online version of Thue's theorem remains
unknown.
\end{abstract}

\section{Introduction}

A \emph{repetition} of \emph{size} $t\geq 1$ in a sequence $%
S=s_{1}s_{2}\ldots s_{n}$ is a subsequence of the form $s_{i+1}s_{i+2}\ldots
s_{i+2t}$ satisfying $s_{i+j}=s_{i+t+j}$ for all $j=1,2,\ldots ,t$. A
sequence $S$ is \emph{nonrepetitive} if there is no repetition (of any size)
in $S$. A celebrated theorem of Thue \cite{Thue}\ (cf. \cite{BerstelThue})
from 1906 asserts that there are arbitrarily long nonrepetitive sequences
over the set of just $3$ symbols. This result has lots of applications and
generalizations (cf. \cite{AlloucheShallit}, \cite{BerstelPerrin}, \cite%
{GrytczukDM}, \cite{Lothaire}, \cite{PezarskiZmarz}). Recently, some game
theoretic variants has been introduced leading to new challenges in the area 
\cite{Pegden}. A basic idea is that a nonrepetitive sequence is created by
two players, Alice and Bob, say, but only Alice cares of avoiding
repetitions. For instance, they may be picking alternately symbols from a
fixed set $A$ and appending them at the end of the existing sequence.
Whenever a repetition occurs, its second part is being erased immediately.
It is proved in \cite{GKM}\ that Alice can still create arbitrarily long
nonrepetitive sequences (no matter what Bob is doing) provided the size of $A
$ is at least $8$.

In this paper we introduce another Thue type game, which we call the \emph{%
online Thue game}. In one round of the game Bob chooses a position in the
existing sequence $S$, which is specified by a number $i\in \{0,1,\ldots ,n\}
$, and Alice is picking a symbol $x\in A$ which is inserted right after $%
s_{i}$, thereby giving a new sequence $S^{\prime }=s_{1}\ldots
s_{i}xs_{i+1}\ldots s_{n}$, with $i=0$ meaning that $x$ is placed at the
beginning of $S$. The goal of Bob is to force Alice to create a repetition,
while Alice will try to avoid that for as long as possible. For instance, if 
$A=\{a,b,c\}$ and $S=acbc$, then Bob catches Alice in one move by choosing $%
i=1$. Indeed, picking any $x\in A$ results in a repetition in $S^{\prime }$: 
$\boldsymbol{aa}cbc$, $a\boldsymbol{bcbc}$, $a\boldsymbol{cc}bc$. Actually,
it is easy to check that Alice has no chance for a longer play over $3$
symbols. A bit more effort is needed to check that the same is true when $A$
is of size $4$, and one may start thinking that there is no finite bound at
all. However we shall prove the following theorem.

\begin{theorem}
There is a strategy guaranteeing Alice arbitrarily long play in the online
Thue game on the set of $12$ symbols.
\end{theorem}

The proof is based on a former result of K\"{u}ndegen and Pelsmajer \cite%
{KundgenPelsmajer}, and Bar\'{a}t and Varj\'{u} \cite{BaratVarju}, on
nonrepetitive coloring of outerplanar graphs. A vertex coloring of a graph $G
$ is \emph{nonrepetitive} if no repetition can be found along any simple
path of $G$.

\begin{theorem}
\emph{(\cite{BaratVarju}, \cite{KundgenPelsmajer}) }Every outerplanar graph
has a nonrepetitive coloring using $12$ colors.
\end{theorem}

The proof of Theorem 1, in a slightly different setting, is given in the
next section. The last section contains discussion and some open problems.

\section{Proof of the main result}

We start with introducing an equivalent setting for the online Thue game
that will be more convenient for our purposes. First notice that this game
can be played on the real line in such a way that Bob is choosing points of
the line and Alice is coloring them using $A$ as the set of colors. Then
after $n$ rounds we have a sequence of points $B=b_{1}b_{2}\ldots b_{n}$, $%
b_{i}\in \mathbb{R}$, written in increasing order $b_{1}<b_{2}<\ldots <b_{n}$%
, and the corresponding sequence of colors $S=c(b_{1})c(b_{2})\ldots c(b_{n})
$, $c(b_{i})\in A$, assigned to the points $b_{i}$ by Alice. The goal of
Alice is to avoid repetitions in $S$.

\begin{proof}
Without loss of generality we may assume that Bob starts with $0$, and
whenever he wants to extend the sequence $B$ to the left or to the right, he
always chooses the closest integer point. So, in the second move he picks
either $-1$ or $1$. In consequence the extreme points $b_{1}$ and $b_{n}$ of 
$B$ are always integers. Furthermore we assume that when Bob decides to
insert a point between $b_{i}$ and $b_{i+1}$, then the new point lies
exactly in the middle of these two. So, if he has chosen $1$ in the second
round, then he may now choose $-1,1/2,$ or $2$. The set of points accessible
in this way is just the set of \emph{dyadic rationals} $D$ which consists of
all numbers of the form $a/2^{k}$, where $a$ is any integer, and $k$ is any
nonnegative integer. We assume that $(a,2^{k})=1$, and we call the exponent $%
k$ the \emph{depth} of $a/2^{k}$. Denote by $D_{k}$ the set of all dyadic
rationals of depth $k$. Notice that $D_{0}$ is just the set\ of all integers.

Now, for each integer $k\geq 0$ consider a graph $P_{k}$ whose vertex set is
the union $D_{0}\cup D_{1}\cup \ldots \cup D_{k}$, with two points joined by
an edge if and only if there is no other point of the union between them on
the line. Clearly $P_{k}$ is a bi-infinite path. Consider finally a graph $G$
on the vertex set $D$ whose set of edges is the union of sets of edges of
all paths $P_{k}$. It is not hard to see that $G$ is an outerplanar graph.
Indeed, we may draw the edges of each path $P_{k}$ as half-circles of
diameter $1/2^{k}$ lying above the line. Then no two edges of $G$ can cross,
and clearly all vertices of $G$ touch the outer face of this embedding.

A simple (but crucial) observation is that any sequence $B=b_{1}b_{2}\ldots
b_{n}$ that can be formed by Bob during the game is a path in $G$. Actually, 
$B$ is a path in a subgraph $G_{n}$ of $G$ which is the union of paths $%
P_{0},P_{1},\ldots ,P_{n}$. Indeed, assume inductively that this holds for $n
$ and consider a sequence $B^{\prime }$ from the next round. Suppose first
that a new point $x$ is placed at the beginning of $B$, that is $B^{\prime
}=xB$. Then, accordingly to the rules, $x$ must be the closest integer to $%
b_{1}$, which itself is an integer. Hence $xb_{1}$ is an edge in $P_{0}$,
and consequently $B^{\prime }$ is a path in $G_{n+1}$. The same holds for $%
B^{\prime }=Bx$. Now, suppose $x$ has been inserted between $b_{i}$ and $%
b_{i+1}$. By inductive assumption these two points lie consecutively on some
path $P_{j}$, $0\leq j\leq n$. Hence, $b_{i},x,b_{i+1}$ lie consecutively on
path $P_{j+1}$. In other words, the edge $b_{i}b_{i+1}$ has been subdivided,
and therefore $B^{\prime }$ is a path in $G_{n+1}$.

To complete the proof we need only to demonstrate that graph $G$ has a
nonrepetitive coloring using $12$ colors. But this follows easily from
Theorem 2 by invoking the compactness principle.
\end{proof}

\section{Remarks and open problems}

First notice that we have actually proved a stronger result asserting that
Alice has a strategy to play \emph{infinitely} long, not just to play \emph{%
arbitrarily} long. The strategy seems nonconstructive, as it is based on
nonrepetitive coloring of an infinite graph obtained via compactness
principle. However, looking into details of the proof of Theorem 2, one has
an impression that it can be adapted to infinite graphs so that an explicit
formula for the nonrepetitive coloring of $G$ is perhaps possible. 

Notice also that in fact we need a weaker coloring in which only
forward-directed paths are nonrepetitive. We suspect that this can be
achieved with smaller number of colors, perhaps $9$ are sufficient.

Another observation leads to a major open problem of the area of
nonrepetitive colorings of graphs. Notice that graph $G$ from the proof of
Theorem 1 has \emph{page number} $1$. What about nonrepetitive coloring of
graphs with higher page number?

\begin{conjecture}
There is a finite constant $N$ such that every graph of page number $2$ has
a nonrepetitive coloring using $N$ colors.
\end{conjecture}

This innocently looking question was propounded some years ago by Idziak.
However, now we know that it is strong enough to imply the following
statement.

\begin{conjecture}
There is a finite constant $M$ such that every planar graph has a
nonrepetitive coloring using $M$ colors.
\end{conjecture}

The best result to date \cite{Dujmovic2} gives a logarithmic upper bound (in
the number of vertices of a graph). Perhaps it would be easier to verify
directed versions of the above problems.

Finally, one may consider online versions of other theorems in combinatorics
on words. Particularly interesting seems the case of \emph{abelian
repetitions}, that is subsequences of the form $r_{1}\ldots r_{t}r_{\sigma
(1)}\ldots r_{\sigma (t)}$, where $\sigma $ is any permutation of the set $%
\{1,2,\ldots ,t\}$. Answering a question of Erd\H{o}s \cite{Erdos}, Ker\"{a}%
nen \cite{Keranen} proved that there exist arbitrarily long sequences over $4
$ symbols avoiding abelian repetitions of any size. Is the online version of
this result possible?

\begin{acknowledgement}
J. Grytczuk acknowledges a support of the Polish Ministry of Science and
Higher Education, grant N206257035.
\end{acknowledgement}


\begin{thebibliography}{99}
\bibitem{AlloucheShallit} J.-P. Allouche, J. Shallit, Automatic Sequences.
Theory, Applications, Generalizations, Cambridge University Press,
Cambridge, 2003.

\bibitem{AGHR} N. Alon, J. Grytczuk, M. Ha\l uszczak, O. Riordan,
Non-repetitive colorings of graphs, Random Structures Algorithms 21 (2002)
336--346.

\bibitem{BaratVarju} J. Bar\'{a}t, P. P. Varj\'{u}. On square-free vertex
colorings of graphs. Studia Sci. Math. Hungar. 44 (2007) 411--422.

\bibitem{BerstelThue} J. Berstel, Axel Thue's papers on repetitions in
words: a translation, Publications du LaCIM, vol. 20, Universit\'{e} du Qu%
\'{e}bec a Montr\'{e}al, 1995.

\bibitem{BerstelPerrin} J. Berstel, D. Perrin, The origins of combinatorics
on words, Europ. J. Combin. 28 (2007) 996--1022.

\bibitem{Dujmovic2} V. Dujmovi\'{c}, F. Frati, G. Joret, D. R. Wood,
Nonrepetitive colourings of planar graphs with $O(\log n)$ colours. Arxiv
manuscript (2012).

\bibitem{Erdos} P. Erd\H{o}s, Some unsolved problems, Magyar Tud. Akad. Mat.
Kutato. Int. Kozl. 6 (1961) 221--254.

\bibitem{GrytczukDM} J. Grytczuk, Thue type problems for graphs, points, and
numbers. Discrete Math. 308 (2008) 4419--4429.

\bibitem{GKM} J. Grytczuk, J. Kozik, P. Micek, New approach to nonrepetitive
sequences. Random Structures Algorithms, DOI 10.1002/rsa.20411.

\bibitem{Keranen} V. Ker\"{a}nen, Abelian squares are avoidable on 4
letters, Automata, Languages and Programming: Lecture Notes in Computer
science, vol. 623, Springer, Berlin, 1992, pp. 41--52.

\bibitem{KundgenPelsmajer} A. K\"{u}ndgen and M. J. Pelsmajer, Nonrepetitive
colorings of graphs of bounded treewidth, Discrete Math. 308 (2008),
4473--4478.

\bibitem{Lothaire} M. Lothaire, Combinatorics on Words, Addison-Wesley,
Reading, MA, 1983.

\bibitem{Pegden} W. Pegden, Highly nonrepetitive sequences: Winning
strategies from the local lemma, Random Structures Algorithms 38 (2011)
140--161.

\bibitem{PezarskiZmarz} A. Pezarski and M. Zmarz, Non-repetitive $3$%
-coloring of subdivided graphs, Electron. J. Combin. 16 (2009) 7.

\bibitem{Thue} A. Thue, \"{U}ber unendliche Zeichenreichen, Norske Vid.
Selsk. Skr., I Mat. Nat. Kl., Christiania 7 (1906) 1--22.
\end{thebibliography}
\end{document}